\documentclass{article}

\usepackage{amssymb,amsmath}
 \usepackage[dvips]{graphicx}
\usepackage{verbatim}

\begin{document}
\begin{center}
{\Large \bf
% Begin Title REQUIRED
Distance Configurations of Points in a Plane with a Galois group
that is not Soluble\\
\medskip

J C Owen and S C Power

% End Title
}\vspace{.5in}

% Begin Name REQUIRED
% J C Owen and S C Power
% End Name

% Begin Address REQUIRED
D-Cubed Limited Park House, \\
Castle Park, Cambridge UK\\
% End Address
\medskip

Department of Mathematics and Statistics,\\ Lancaster University,
Lancaster, UK
% Begin E-Mail OPTIONAL

% End E-Mail
\end{center}

%J C Owen D-Cubed Limited Park House, Castle Park, Cambridge UK

%S C Power Department of Mathematics, Lancaster University,
%Lancaster, UK

\begin{abstract}
We have conjectured that the constraint equations defined by a generic
Laman graph are not soluble by radicals when the graph is $3$-connected.
We prove that this conjecture follows from the following simpler
conjecture: the constraint equations defined by a generic Laman graph
are not soluble by radicals if the graph does not contain a proper
subgraph which is itself a Laman graph.
\end{abstract}

\section{Introduction}

 We consider the algebraic equations that result from specifying
the relative distances between points in a plane. This gives a
system of quadratic equations of the form
\[
(x_i-x_j)^2+(y_i-y_j)^2 = d_{ij}, \quad \quad \quad \quad \quad (1)
\]
where there is a pair of variables $(x_i,y_i)$ corresponding to
each point $i$ and an equation for each of the specified
(squared) distances
$d_{ij}$
 between points $i$ and $j$. The distances are taken as
input parameters.

This is a simplified model of the more general problem that occurs
in Mechanical Computer Aided Design (MCAD), namely the
determination of the location of various types of geometric
elements (points, line segments and arcs) when the designer
specifies the relative distances between them. A similar system of
equations arises in the study of rigid structures  \cite{whi}, in
molecular structure \cite{hav-cri} and in robotics \cite{str}.
The configurations
that are of special interest are what engineering designers call
properly dimensioned or minimally rigid which means that the
number of distances is just sufficient to ensure discrete
solutions for the points (subject to the obvious rigid body
motions in the Euclidean plane). The general problem of geometric
constraint solving is of considerable practical interest and over
the past decade almost all MCAD products have included some form
of geometric constraint solving. The majority of them use
algorithms which are based on the mathematics described here \cite{dcm}.

The sets of equations which are best suited to interactive or
real-time solution are those that can be solved without recourse
to numerical iteration methods. Examples include a system of
linear equations (which can be solved using Gaussian elimination)
or a system of non-linear polynomial equations which can be
transformed into triangular form involving at most the solutions
of quadratic equations (quadratically soluble, QS) or equations
which can be solved by extracting higher roots (radically soluble,
RS). As well as being of practical importance, this classification
leads to a nice connection with classical mathematics. QS
equations correspond to geometric configurations that can be
constructed using only a ruler and compass and the zeros of RS
equations define an extension field with a Galois group that is a
soluble group.

The equations (1) above are invariant under the continuous
group of rotations and translations. This invariance can be removed
by specifying as additional input parameters
the coordinates of any pair of points which are joined by
a distance \cite{owe-pow}.
All such systems are equivalent by an affine transformation
(over input parameters)
and so the Galois group of a solution to the equations does not depend
upon the pair of points which are chosen. We will simply assume
that some convenient pair have been chosen.

The determination of the Galois group for an arbitrary geometric
constraint problem seems to be very difficult. There are certainly
many quite apparently complex configurations that are QS \cite{dcm} and
there are also some quite simple ones (involving for example just
six points in a plane) whose Galois group is not soluble and whose
equations are neither QS nor RS \cite{owe-pow} and Section 3 below.

Although any instance of constraint solving will
typically assign a set of
rational values to the distances, we are commonly interested in
finding solution methods which are independent of the particular
values used for the distances. This leads us to consider solution
methods when the distances are assigned values which are
algebraically independent (over the rationals) which we will also
call generic. We have shown in \cite{owe-pow} that the form of equations (1)
means that for generic distances all of the
(finitely many) solutions
have the same Galois group.

An abstract graph can be defined to represent a set of points with
generic relative distances by corresponding the points with
vertices of the graph and the generic distances with edges of the
graph. We consider a configuration of $N$ points in a plane with
$E$ relative distances. The corresponding abstract graph has $N$ vertices
and $E$ edges. The graph is described as independent if every
subgraph with $n$ vertices and e edges satisfies $2n-e\ge 3$ and
maximally independent if additionally we have $2N-E=3$. It is a
well-known result of Laman \cite{lam} that maximal independence for a
graph is equivalent to infinitesimal rigidity for the
corresponding generic framework and that these are also equivalent
to the ridigity of the framework as a bar-joint mechanism \cite{asi-rot}. For
this reason, a graph which is maximally independent is also known
as a Laman graph and we will sometimes adopt this shorthand
terminology. We note that in the usual definition of a generic framework
the coordinates of the points (or joints) are taken as the
algebraic independents. However, we have shown in \cite{owe-pow} that this is
equivalent to considering the distances as algebraic independents.
In view of this equivalence we may simply refer to a Laman graph
as generic if it represents a system of constraint equations with
generic distances and we may also refer to the Galois group of a
generic Laman graph as a shorthand for the Galois group of the
corresponding equations with generic distances.

Some years ago one of us proposed that for generic distances there
is also an abstract graph characterisation of whether a
configuration has solutions that are QS \cite{owe}. We showed that any
Laman graph has a recursive decomposition into subgraphs
which are also maximally independent (after the possible
introduction of some virtual edges) and are either 3-cycles (also
called triangles) or are vertex $3$-connected. We showed that the
graph is QS (and thus RS) if these decomposition subgraphs are
all 3-cycles, and proposed that otherwise the graph is not-QS. We
have extended this to conjecture that the graph is also not-RS
\cite{owe-pow}.
In other words, we have proved a sufficient condition for any
generic graph to be RS and we propose that this condition is also
sufficient. This sufficient condition would follow from the truth
of the following conjecture:
\medskip

\noindent {\bf Conjecture A:} {\it  If a generic Laman graph is $3$-connected
then it is not solvable by radicals (not-RS)}.
\medskip

We still do not have a complete proof of this conjecture. We have
proved the conjecture for the (infinite) class of graphs that have
a planar embedding (that is, the graph can be drawn on a plane
with no edge crossings) \cite{owe-pow}. Although our initial focus on planar
graphs was made in order to identify a class of graphs for which
we could complete a proof, we have subsequently found that planar
Laman graphs are of specific interest in robotics problems \cite{str} and
so even this partial proof has some merit. Nevertheless the
introduction of planar graphs does not seem to be intrinsic to the
problem.

The main new result which we present in this paper is that the
Conjecture A follows from a simpler conjecture which refers to a
smaller class of graphs, namely:
\medskip

\noindent {\bf Conjecture B:}  {\it  If a generic Laman graph does not
contain a proper maximally independent subgraph then it is not
solvable by radicals (non-RS)}.
\medskip

In order to appreciate the significance of Conjecture B we
introduce the notion of a basic Laman graph as follows:

\medskip
\noindent {\bf Definition:}{\it  A basic Laman graph has $2N-E=3$ with a
corresponding strict inequality $2n-e>3$ in every proper subgraph
with $n>2$.}
\medskip

%SMOOTH NEXT PARAG. NOT CLEAR - "CIRCUIT" IN A "LAMAN PLUS ONE"
%GRAPH SEEMS TO NEED TWO DEFINITIONS
%We note the similarity between our definition of a basic Laman
%graph and the definition of a circuit in a Laman-plus-one graph
%[9] which requires $2N-E=2$ and $2n-e>2$ in every proper subgraph.
%An edge can be added between any pair of vertices in a basic Laman
%graph to give a Laman-plus-one circuit. Recently a constructive
%method has been given for generating all Laman-plus-one circuits
%starting from $K_4$ (the complete graph on 4 vertices) [10]. It is
%interesting to speculate whether there is a corresponding
%construction for all basic Laman graphs.

Conjecture B can now be restated in the intuitively appealing
form:
\medskip

\noindent {\bf Conjecture B:} {\it Basic Laman graphs are not generically
solvable by radicals (not-RS)}.
\medskip

It is easy to show (see corollary to Lemma 1 below) that every
basic Laman graph is $3$-connected and non-planar. The converse is
not true. This is why Conjecture B (even when it is taken to refer
only to non-planar graphs) refers to a smaller class of graphs
than Conjecture A. There is one basic Laman graph with $N=6$,
there are none with $N=7$ and two with $N=8$. These are shown in
Figure 1. Figure 2 shows a non-planar graph with $N=7$ which is
$3$-connected but not basic.

Section 2 of the paper is used to prove Lemma 1 and the following
main theorem
\medskip

\noindent {\bf Theorem 1}: {\it Conjecture B implies Conjecture A.}
\medskip

In section 3 we show by explicit calculation using Maple that the
smallest basic Laman graph is not-RS. This adds weight to our
conjecture that they are in fact all not-RS.

\section{Proof of The Main Theorem}

We use the same definitions as in \cite{owe-pow} for the concepts which we
adopt from graph theory and algebraic geometry.  The reader is
referred to that paper for complete definitions. For the sake of
readability we reproduce our main graph theory definitions here.

A graph $G$ is a set of vertices $x$ and a set of edges $(xy)$
such that the edge $(xy)$ is in $G$ only if the vertices $x$ and
$y $ are in $G$. The {\it order} $|G|$ is the number of vertices in $G$.
A vertex $x$ is described as {\it incident}  to the edge $(xy)$. Vertices
$x$ and $y$ are {\it adjacent} in $G$ if $G$ contains the edge $(xy)$.
The {\it union} and {\it intersection} of graphs $G$ and $H$ is the union and
intersection of the sets (with duplicates removed). A subgraph
$S$ of $G$ is described as {\it proper} if it is contained in $G$ and
not equal to $G$ and has at least three vertices. A subgraph $S$
of $G$ is {\it vertex induced} if, whenever $S$ contains $x$ and $y$ and
$G$ contains $(xy)$ then $S$ also contains $(xy)$. The graph
$G\backslash S$ is the graph induced by the vertices of $G$ which
are not in $S$. The {\it vertices of attachment} of $S$ in $G$ are the
vertices of $S$ which are adjacent in $G$ to vertices in
$G\backslash S.$ An {\it internal vertex} of $S$ is a vertex
of $S$ which is
not a vertex of attachment. A sub-set $V$ of vertices of $G$ is
said to {\it separate} $G$ if there exists a pair of vertices $\{a, b\}$
in $G\backslash  V$ such that every connected path of edges which
joins $a$ and $b$ always includes a member of $V$. A graph is
{\it m-connected}
if $|G|>m$ and $G$ has no separation sets of order $m-1$. The disconnected
components $\{C_i\}$ of the graph induced
by the vertices of $G\backslash V$ are the
separation components of $G$ with respect to $V$ \cite{die}. The
subgraphs of $G$
induced by the vertices of $C_i$ and $V$ are the {\it separation blocks}
 $\{S_i\}$ of
$G$ with
respect to $V$. Thus if a graph is $m$-connected, but not $m+1$-connected,
there is a set of at least two separation blocks $\{S_i\}$ such that $|S_i|>m$,
the vertices of $S_i\cap S_j$ are exactly the vertices of $V$ for any $i$ and
$j$ and
$G = \cup_i S_i$.

If a graph $G$ has $n$ vertices and $e$
edges then the {\it freedom number}, $free(G)$  is  $ 2n-e-3$.

A graph $G$ is {\it independent} if $G$ and every proper subgraph of
$G$ have a non-negative freedom number. $G$ is a {\it maximally
independent} graph or equivalently a {\it Laman} graph if in addition
$free(G)=0$. A Laman graph is {\it basic} if it has no proper subgraph
which is maximally independent.

The properties of a basic Laman graph which we have described
above are easily obtained from the following:
\medskip

\noindent {\bf Lemma 1}: {\it  Let $G$ be a Laman graph with a vertex separation
pair $\{a,b\}$ which separates $G$ into at least two
separation blocks $S_i$ such that $G = \cup_i S_i$ and
the subgraphs
$S_i\cap S_j$
have vertex set $\{a,b\}$.
% and edge set $\{(ab)\}$ if and only
%if $G$ contains the edge $(ab)$.
If $G$ contains
$(ab)$ then $free(S_i) = 0$, for all $i$, otherwise there is a block
$S_1$ such that $free(S_1) = 0$ and $free(S_j)=1$ all $j>1$.}
\medskip

\noindent {\it Proof.}
 Suppose that
disjoint graphs $H$ and $K$ are joined to create $H \cup K$ by
identifying a pair of vertices $\{a,b\}$ from each and removing
repetition of edge $(ab)$ if both contain $(ab)$. If neither graph
has edge $(ab)$ then $free(H \cup K) = free(H) +
free(K) -1.$ It follows (by repeated joining) that if $G$ does not
contain $(ab)$ then $free(G)$ is the sum of the freedom numbers of
the components minus $k-1$ where $k$ is the number of components.
Since the components are independent and so have non-negative
freedom number, the second half of the Lemma follows. When $G$
contains $(ab)$ there is a similar argument.

\hfill $\Box$
\medskip

\noindent {\bf Corollary}: {\it Every basic Laman graph $G$ with $|G|>3$ is $3$-connected
and non-planar.}
\medskip

\noindent {\it Proof.}  Every Laman graph $G$ is 2-connected
%[\cite{owe}, Lemma 4.8] ??
so if
$G$ is not $3$-connected with $|G|>3$ it has a separation pair and thus
a proper subgraph $S_1$ with $free(S_1) = 0$. Every planar Laman graph
$G$ contains a 3-cycle $C$, by Corollary 4.10 of \cite{owe-pow}, for which $free(C)
= 0$, and $C$ is a proper subgraph for $|G|>3$.
\hfill $\Box$

\medskip

Our strategy to prove Theorem 1 is to show that if there is a
$3$-connected, non-basic Laman graph with $N>6$ vertices which is
radically soluble (RS), then there is a $3$-connected Laman graph
with fewer than $N$ vertices which is also RS. This defines a
sequence of $3$-connected Laman graphs such that each successor
graph is RS if its ancestor is RS. The sequence terminates on
either a non-basic Laman graph with $N=6 $ or on a basic Laman graph.
The only non-basic Laman graph with $N=6$ is the doublet, which we
have proved to be not-RS in [5] and so Conjecture A is proved true
if Conjecture B is true.

The main graph reduction step which we use is edge contraction
\cite{die}. This is illustrated in figure 3. If $G$ is a graph with an
edge $(xy)$ which joins vertices $x$ and $y$ then the edge contracted
graph $G/(xy)$ is obtained from $G$ by deleting the edge $(xy)$,
identifying the vertices $x$ and $y$ and deleting any duplicate edges
which may result. The fundamental result which we need from
algebraic geometry is proved in Theorem 6.1 of \cite{owe-pow}, namely:
\medskip

{\it If both $G$ and $G/e$ are maximally independent graphs then $G/e$ is
generically soluble by radicals (RS) if $G$ is generically soluble
by radicals (RS).}
\medskip

In order to use this graph reduction we must find edge
contractions which maintain both maximal independence and
$3$-connectivity. It is easy to see that a necessary condition to
preserve maximal independence when an edge in a graph $G$ is
contracted is that the edge should be in exactly one 3-cycle of $G$.
In \cite{owe-pow} we focused on planar graphs in order to ensure that $G$
contains a 3-cycle. We say that an edge $e$ of a maximally
independent graph $G$ is contractible if $G/e$ is maximally
independent.

To prove Theorem 1 when $G$ may not contain a 3-cycle we first
perform some surgery on $G$ in order to generate a graph which
does have a 3-cycle. Since $G$ can be assumed to have a subgraph
$R$ which is maximally independent we can replace this subgraph in $G$
with a subgraph which is a simple triangulation of the vertices
of attachment of $R$ in $G$. This surgery is illustrated in Figure 4
and defined precisely in Lemma 4 below.

This subgraph replacement is a special case of a more general
replacement scheme in which a maximally independent subgraph $R$ is
replaced with another maximally independent graph $R'$ using the
same vertices of attachment. This preserves many important
properties of $G$ as shown by the following.
\medskip

\noindent {\bf Lemma 2}. {\it Let $G$ be a graph
with a maximally independent proper
vertex induced
subgraph $R$ which has vertices of attachment $\{c_1,...,c_m\}$ in $G$. Let
$R'$ be any maximally independent graph with $n \ge m$ vertices such
that m of its vertices are identified with the vertices $\{c_1,\dots ,c_m\}$ in
$G$ and the remaining vertices are called internal vertices of $R'$. A
graph $G'$ is obtained from $G$ by deleting all of the edges and all
of the internal vertices of $R$ and replacing them with the edges
and internal vertices of $R'$. Then (i) $free(G')=free(G)$. (ii)    If
$G$ is independent then $G'$ is independent. (iii)   If $G$ is
maximally independent and both $G$ and $R'$ are generically RS, then
$G'$ is generically RS.}
\medskip

\noindent {\it Proof.} Suppose that $G, G', R$ and $R'$
have $N, N', n$ and $n'$ vertices
and $E, E', e$ and $e'$ edges. Then $N' = N-n+n', E'=E-e+e'$ and
$free(G')=2N'-E'-3 = 2N-E + 2n'-e'-(2n-e)-3 = 2N-E-3=free(G)$.

To prove (ii) we show that every subgraph $S'$ of $G'$ has
$free(S')\ge 0$.

Suppose to the contrary that there is a subgraph $S'$ of $G'$ with
$free(S')<0.$ Then $S'$ contains at least $2$ of the vertices
$c_i$ or else it would also be a subgraph of $G$.  Then $free(S'
\cup R') = free(S')+free(R') - free(S'\cap R')$ and $free(S'\cap
R')\ge 0$ because $(S'\cap R')$ is a subgraph of $R$, $free(R')=0$
and $free(S')$ is negative. This implies $free(S' \cup R')<0$.

$(S' \cup  R')$ is a graph which contains the maximally independent
subgraph $R'$. If $R'$ is replaced by the maximally independent graph
$R$ in the manner described in the Lemma statement above then the
resulting graph has the form $(S \cup R)$ for some subgraph $S$ of $G$
and by the first assertion of the lemma $free(S \cup R)= free(S' \cup
R')<0$. This contradicts the requirement that $G$ is independent
since $(S \cup R)$ is a subgraph of $G$.

To prove (iii) we appeal to some standard results of algebraic
geometry \cite{cox-lit-osh}.

Let the (squared) distances in $R$ be $\{d_r\}$, the distances in $R'$ be
$\{d_s)$, distances in $G\backslash R$ and in $G\backslash
R'$ be $\{d_t\}$, and the distances in $G$ be $\{d_g\}$.
Each set of generic distances is thus an
algebraically independent set of indeterminates and in particular
the (transcendental) field extension $\mathbb{Q}(\{d_g\})$ is a
base field over which we can form the algebraic closure.
It follows from the maximal independence
of $R$ and $G$ that their equation sets determine algebraic
varieties $V(R)$ and $V(G)$ (over the appropriate algebraic
closures) which consist of finitely many points ([4], [16]), and
so (by definition) are zero dimensional varieties.

\begin{comment}
The equations of $R$ generate a polynomial ideal, $I(R)$ say, in
the algebra of polynomials over $\mathbb{Q}$ in the coordinate
variables and the distances $\{d_g\}$. This ideal is certainly
contained in the elimination ideal $I(G)_R$ which, by definition,
is the set of polynomials in $I(G)$ which involve only the
coordinate and distance variables of $R$.
 Both $I(R)$ and $I(G)_R$
have varieties with a finite number of zeros
 because the varieties of $R$ and $G$
both have dimension $0$.

For any ideal $J$ let $V(J)$ be the variety of zeros of $J$ in some
suitable algebraic closure (or in the "universal" field). Since
$I(R)$ is contained in any $I(G)$, any zero of $I(G)$ also gives a zero
of $I(R)$ simply by restricting to the variables of $R$.
It will be convenient to use the  notation
$I(G|w)$ to indicate the ideal generated by
the equations of $G$ evaluated at the point $w$.
\end{comment}

The equation set $R$ is contained in the equation set $G$ and so any point
of $V(G)$ gives a point of $V(R)$ simply by restricting to the variables $R$.
It will be convenient to use the notation $G|w$ to indicate the equation
set $G $ with (some of) its variables evaluated at the point $w$ and
$\pi_{G\backslash R}(V(G))$
to represent the points of the variety $V(G)$ projected onto the subset of
points represented by the variables in the equation set $G\backslash R.$

We have
%\pi(V(G)) = Uw_in_V(R)(V(G|w))=20
%            =3D Uw_in_V(R)(V(G\R)|w).
\begin{eqnarray*}
\pi_{G\backslash R}(V(G)) & = \cup_{w \in V(R)}
(V(G|w))\\
& = \cup_{w \in V(R)}
(V((G\backslash R)|w)).
\end{eqnarray*}
Now $V(G)$ is radical over $\mathbb{Q}(\{d_r\},\{d_t\})$
and so $\pi_{G\backslash R}(V(G))$
 is radical over $\mathbb{Q}(\{d_r\},\{d_t\})$ and for any
 zero $w$ in $V(R)$
we have that $V((G\backslash R)|w)
$ is
radical over $\mathbb{Q}(\{d_r\},\{d_t\})$. Let $(x_1,...,x_{2k})$
($k \ge m)$ be a zero of $V(R)$ where $x_{2i},x_{2i+1}$ are the coordinates
of $c_i$.
Then
$\mathbb{Q}(\{d_r\})$ is contained in $\mathbb{Q}(\{x_1,...,x_{2k}\})$ so
$V((G\backslash R)|w)$
is radical over $\mathbb{Q}(\{x_1,...,x_{2k}\},\{d_t\})$. But
$\{x_1,...,x_{2k}\}$ are
algebraically independent over $\mathbb{Q}$ and  so
$V((G\backslash
R)|(x_1,...,x_{2k}))$
is radical over \\
$\mathbb{Q}(\{x_1,...,x_{2k}\},\{d_t\})$.
This holds for any algebraically independent set
$\{x_1,...,x_{2k}\}$.

\begin{comment}
*************************************************

In fact these ideas are equal because $I(R)$ is prime over
$\mathbb{Q}(\{d_r\})$ and
thus prime over $\mathbb{Q}(\{d_g\})$. $I(R)$
is therefore a maximal ideal in $\mathbb{Q}(\{d_g\})$ and
hence equal to $I(G)_R$ since it is contained in $I(G)_R$.

For any ideal $J$ let $V(J)$ be the variety of zeros of $J$ in some suitable
algebraic closure (or in the "universal" field). Then by the extension
theorem

$\pi_{G\backslash R}(V(I(G)) = \cup_{\mbox{zeros of }~ I(G)_R}
\{V(I(G))\}|\mbox{zeros of }~ I(G)_R$

$                 = \cup_{\mbox{zeros of }~
I(R)}\{V(I(G))\}|\mbox{zeros of}~ I(R)$

$                 = \cup_{\mbox{zeros of }~
I(R)}\{V(I(G\backslash R))\}|\mbox{zeros of}~ I(R)$

BELOW TO LATEX PROPERLY
\end{comment}

To construct a radical zero of $G'$, we select the vertices $c_1$ and
$c_2$ as the base vertices \cite{owe-pow} and obtain the coordinates for the
vertices $c_3,\dots ,c_m$ by solving for the coordinates of the vertices $c_3$,
$c_4$ etc in sequential order where the coordinates for the vertex $c_j$
is determined from the equations represented by the edges $(c_jc_{j-1})$
and $(c_jc_1)$. This can be done for each vertex in turn by solving a
single quadratic equation \cite{owe}. Thus each of the coordinates of
${c_1,\dots ,c_m}$ are in a radical extension of  $\mathbb{Q}(\{d_s'\})$ and these
coordinates are algebraically independent because $R'$ is also a
generic Laman graph. The coordinates of the remaining points of $G'$
are determined as zeros of the equations of $G'\backslash R'$ which are the
same as the zeros of the equations of $G\backslash R$ and therefore lie in a
radical extension of $\mathbb{Q}(\{x_1,...,x_{2m})\},\{d_t\})$ which is in turn in a
radical extension of $\mathbb{Q}(\{d_s'\},\{d_t\})$.
\hfill $\Box$

\medskip

If this surgery is performed on an arbitrary Laman subgraph then
the resulting graph may not be $3$-connected. This is
illustrated in Figure 5, where the graph $G$ is $3$-connected
but the graph $ G\# R$ is not. Notice however that the subgraph
$R$ can be expanded to include the three vertices at the bottom
 of the graph. The resulting subgraph $T$ is maximally independent
  and has 3 vertices of connection in $G$. The graph $G\# T$ is the doublet
  which is $3$-connected.  These considerations lead us to restrict
  attention to subgraphs which are maximal according to the
  following definition:

A proper subgraph $R$ of a graph $G$ is maximal with respect to
property $P$ if there is no proper subgraph $S$ of $G$ with the
property $P$ such that $R$ is a proper subgraph of S.

In order to use this property we need the following lemma:
\medskip

\noindent {\bf Lemma 3}: {\it Let $G'$ be the graph which
is obtained from a maximally
independent graph $G$ by the replacement of a maximally
independent proper subgraph $R$ with a maximally independent graph
$R'$, in the manner described in Lemma 2. If $R$ is maximal in $G$
(with respect to the property of maximal independence) then $R'$ is
maximal in $G'$.}

\medskip

\noindent {\it Proof.} Suppose to the contrary that there is a proper subgraph $S'$
of $G'$ which is maximally independent and which contains $R'$ as a
proper subgraph. This means there is a vertex w in $G' \backslash S'$, so w
is also in $G\backslash R$. The vertices and edges of $S'\backslash R'$
are vertices and
edges of $G$ and so replacing the proper subgraph $R'$ in $S'$ by the
graph $R$ gives a graph $S$ which is a proper subgraph of $G$
(proper, because it does not contain w) and which contains $R$ as a
proper subgraph. $S$ is maximally independent by Lemma 2 which
contradicts the fact that $R$ is maximal.
\hfill $\Box$

\medskip

With these preliminaries we can now prove the following useful
properties of the proposed graph surgery.
\medskip

\noindent {\bf Lemma 4}: {\it Let $G$ be a $3$-connected maximally independent graph which has
a maximally independent proper subgraph $R$ with $m>2$ vertices of attachment
 ${c_1,\dots ,c_m}$. Define  $G\# R$
  to be the graph obtained from $G$ by deleting all
 the edges of $R$ and all of the internal vertices of  $R$, adding the cycle
 of edges $(c_1c_2), (c_2c_3),..., (c_mc_1)$ and, for $m>3$, the edges $(c_1c_3),
 (c_1c_4),...,(c_1c_{m-3}), (c_1c_{m-2})$. Then

(i)  $G\# R$ is maximally independent.

If  $R$ is maximal (with respect to the property of maximal
independence) then

(ii)  $G\# R$ is $3$-connected.

(iii) Each of the edges $(c_ic_{i+1})$ is contractible.

(iv) Any proper maximally independent subgraph of  $G\# R$ with an
internal vertex is also a proper maximally independent subgraph
of  $G$ with an internal vertex.
}
\medskip

\noindent {\it Proof.} Point (i) follows immediately from Lemma 2.

Let $R'$ be the subgraph of  $G\# R$ induced by the vertices
${c_1,...,c_m}$. $R'$ is
maximally independent by construction and is maximal in  $G\# R$ by Lemma 3.

To prove (ii), (iii) and (iv) we first show that if $S'$ is a proper maximally
independent subgraph of  $G\# R$ not contained in $R'$,
then $|R'\cap S'| <2$. For if
$|R'\cap S'| >1 $ then by Lemma 4.2 of \cite{owe-pow}
$R' \cup  S'$ is maximally independent and
thus $G$ = $R' \cup  S'$ because $R'$ is maximal. Each vertex $c_i$
is incident to an
edge which is not in $R'$ so this edge and thus each $c_i$ is in $S'$
which implies
that $S'$ contains $R'$ and thus $S'=G'$
which contradicts the requirement that $S'$
is proper.

If  $G\# R$ is not $3$-connected then it has a separation pair $\{a ,b\}$
and both $a$
and $b$ are in $R'$ or else $\{a, b\}$ would also be a separation pair for $G$.
By
Lemma 1 there is a proper subgraph $S_1$ with $free(S_1)=0$ which has exactly
$a$ and $b$ as vertices of attachment in $G$. Thus
$|(R'\cap S_1)|>1$ and $S_1$ is
contained in $R'$ which is impossible because every vertex of $R'$ is incident
to a vertex of $G\backslash R$ and would be a vertex of attachment of $S_1$.

The edge $(c_ic_{i+1})$ is in the 3-cycle $\{c_i,c_{i+1},c_1\}$ of $R'$. By the
point proved above there is no maximally independent subgraph
which is not contained in $R'$ and which also contains $(c_ic_{i+1})$. All
maximally independent subgraphs in $R'$ contain $c_1$.  Thus by Lemma
4.5 of [5] the edge $(c_ic_{i+1})$ is contractible.

If $S'$ is a proper maximally independent
subgraph of  $G\# R$
with an internal vertex then $S'$ is not contained in $R'$ so
$R'\cap S'$ has a most
one vertex and no edges. Thus $S'$ is also a proper subgraph of $G$.
\hfill $\Box$

In order to complete the reduction required to prove Theorem 1, we
have to allow for the fact than even though  $G\# R$ has contractible
edges, the result of a contraction may not be $3$-connected. We have
shown in Lemma 1 that if a Laman graph $G$ has a vertex separation
pair $\{a, b\}$, then $G$ is the union of subgraphs $S_i$, all
with $free(S_i) = 0$ or $1$,
and that at least one of these, say $S_1$,  has $free(S_1)=0$.
If a new (virtual)
edge $(ab)$ is added to every $S_i$ for which $S_i=1$ then each
graph
$S_i'$, defined as
$S_i$ with $(ab)$ added if
$free(S_i) = 1$, is a Laman graph.
\medskip

We have shown in \cite{owe-pow}  that the graph $G$ is
generically RS if and only if each of the $S_i'$ is generically RS.
This
observation forms the basis of our recursive method for generating the
zeros of the generic varieties   of Laman graphs which are RS.

In order to use this observation in the proof of Theorem 1, we
need to carry the analysis a little further because the components
$S_i'$ may not be $3$-connected. This could happen either because $G$
has more that one vertex separation pair or because some subgraph
such as $S_1$ which does not have a virtual edge added to it may have
a separation pair which is not a separation pair of $G$. The former
problem is managed by separating $G$ at all of its separation
pairs. The latter is managed by adding a new virtual edge $(ab)$ to
every subgraph $S_i$ which does not already contain this edge and
defining $B_i(ab) = S_i \cup (ab)$, for all $i$. A new virtual edge $(ab)$ in
$B_i(ab)$ is marked as "redundant" if $free(S_i)=0$. The redundant edges
are included when we consider the $3$-connectivity of $B_i$ but are
excluded when we consider whether $B_i$ is RS. Each $B_i$ is 2-connected
because every Laman graph is 2-connected \cite{owe-pow} . It is easy to see
that any separation pair of any of the $B_i (ab)$ is distinct from
$\{a,b\}$ and that it is also a separation pair of $G$ \cite{hop-tar} . The
converse is also true because $G$ is a Laman graph. This is shown
in the following lemma and leads to a unique decomposition for $G$
which is a little simpler than that given in \cite{hop-tar}  for a general
graph.
\medskip

\noindent {\bf Lemma 5:} {\it  Let $G$ be a Laman graph. Then

(i) If  $\{a, b\}$ and
$\{c, d\}$
are any distinct separation pairs of $G$ which separate $G$
respectively into blocks $B_i(a,b)$ and $B_i(c, d)$ as described
above then $\{c,d\}$ is also a separation pair of some block
$B_1(a,b)$. If $B_1(c,d)$ is the corresponding block of the
separation pair $\{c,d\}$ which contains $\{a,b\}$, then separating $G$ at
$\{a,b\}$ and $\{c,d\}$ gives blocks $B_i(a,b),i>1, B_j(c,d),j>1,$ and
$B_1(a,b)\cap  B_1(c,d)$. These blocks are independent of the
order of the separations.

 (ii)    $G$ has a unique decomposition
into blocks $B_i$ where each block is a subgraph of $G$ plus
virtual edges, some of which are marked as redundant. Each $B_i$ is
either (a) a 3-cycle, or (b) a 3-connected graph .

(iii)   Each $B_i$
(ignoring redundant virtual edges) is maximally independent and
$G$ is generically RS if and only if each $B_i$ is generically RS.

(iv) At least one of the blocks $B_i$ contains no redundant
virtual edges.}
\medskip

\noindent {\it Proof.} If a separation pair $\{c, d\}$ of $G$
is not a separation pair
of any $B_i(a, b)$ then the vertices $c$ and $d$ are in different
blocks $B_i(a,b)$ and $B_j(a,b)$. Then $G$ is the union of four
proper subgraphs as follows. $G = S_1 \cup S_2 \cup S_3 \cup S_4$,
where
$S_1\cap S_2=a, S_2\cap S_3=c,
S_3\cap  S_4=b, S_4\cap  S_1=d $ (see \cite{die}, this
is the case that produces cycles of length greater than $3$ of  original and
virtual edges in the decomposition described in this reference).
Then
$free(G)  = free(S_1)+free(S_2)+free(S_3)+free(S_4)+1$ and this
contradicts $free(G)=0$ because each $free(S_i)\ge 0$. Thus the pair
$\{c,d\}$ is
in some block $B_1(a,b)$ and similarly the pair ${a, b}$
is in some block $B_1(c,d)$. The remaining blocks $B_i(c,d),
i>1$ are all contained in $B_1(a,b)$ because they contain the pair
$\{c,d\}$ and do not contain the pair $\{a,b\}$. Thus the pair $\{c,d\}$
separates $B_1(a,b)$ into $B_i(c,d)$, all $i>1$ and $B_1(a,b)\cap  B_1(c,d).$
Similarly the pair $\{a, b\}$ separates $B_1(c,d)$ into $B_j(a,b)$, all $j>1$
and $B_1(a,b)\cap  B_1(c,d)$. The two pairs $\{a,b\}$ and $\{c,d\}$ separate
$G$ into $B_i(a,b),i>1, B_j(c,d),j>1£$ and $B_1(a,b)\cap  B_1(c,d)$ and this
is independent of the order in which $G$ is separated.

Suppose that $G$ has $m$ distinct separation pairs. The
decomposition in (ii) is obtained by successively separating $G$
and the resulting separation blocks and adding virtual edges
and marking some of them as redundant as described above. This
procedure terminates after $m$ separations when none of the
blocks have any separation pairs. The resulting blocks are
independent of the order of the separations by (i). Each of the
blocks $B_i$ is either $3$-connected or a $3$-cycle of edges since
these are the only 2-connected graphs without a separation pair.

Each block $B_i$ is maximally independent (when redundant virtual
edges are ignored) since this property is maintained for every
block at every separation. Similarly $G$ is RS if and only if
every $B_i$ is RS because this property is also maintained at each
separation.

At most one virtual edge is marked as redundant for each
separation pair of $G$, so there are at most $m$ redundant virtual
edges. However after $m$ separations there are at least $m+1$
blocks and so at least one block has no redundant virtual
edges.
\hfill $\Box$

\medskip

With these preliminary lemmas complete we now give the main graph
reduction lemma which allows us to maintain $3$-connectivity
following edge contractions.
\medskip

\noindent {\bf Lemma 6:} {\it Let $G$
be a $3$-connected, maximally independent graph
$|G|>6$. Suppose that $G$ has no proper maximally independent
subgraph with an internal vertex and that $G$ has a contractible
edge $e$. Then either

(i) $G$ has a contractible edge $f$ such the $G/f$ is $3$-connected, or
(ii)    All separation blocks of $G/e$ (as defined above) are
$3$-connected.}
\medskip

\noindent {\it Proof.} If $e$ is on a $3$-cycle of $G$ which also has a vertex of
degree $3$ then by Lemma 4.7 of \cite{owe-pow}  at least one edge of the $3$-cycle
is contractible to give a $3$-connected graph and (i) is satisfied.
Otherwise we may assume that G/e is not $3$-connected and that the
edge e is on a $3$-cycle all of whose vertices have degree at least
4. Since every vertex of $G$ has degree at least $3$ this implies
that every vertex of $G/e$ has degree at least $3$.

Let $e=(xy)$ . Since $G/e$ is not $3$-connected it has a separation pair
$(v, w)$ and one of the vertices, $v$ say, is the combination of the
vertices $x$ and $y$ from $G$.  This means that the vertex triple
$(x,y,w)$ separates $G$.

We claim that in fact $(x,y,w)$ divides $G/e$ into exactly two
components and that $G$ does not contain the edge $(xw)$ or the edge
$(yw)$. In \cite{owe-pow}  we proved this property (Lemma 4.8) when $G$ is a
planar graph using Kuratowski's theorem. In fact this property is
also true for non-planar graphs which do not contain a proper maximally
independent subgraph. For suppose $(x,y,w)$ divides $G$
into $m$ blocks $S_i$ each with freedom number $f_i$ and that
$d = 0,1$
or $2$ according to how many of the edges $(xw)$ or $(yw)$ are in $G$. The
freedom number of the subgraph induced by the vertices $(x,y,w)$ is
$(2-d)$.
 Then $0 = free(G) = \Sigma_i(f_i)-(m-1)(2-d)$ and so
$\Sigma_i(f_i-2+d) =d-2.$ Each $f_i > 0$ since otherwise the subgraph $S_i$
would be maximally independent with $3$ vertices of attachment in
$G$ and an internal vertex. Thus $\Sigma_i(d-1)\le d-2$ and this
requires $d=0$ and then at least two of the $f_i =1$, say $f_1=f_2=1$. The
subgraph $S_1 \cup  S_2$ has  freedom number $f_1+f_2-2=0$ and is a proper
maximally independent subgraph of $G$ with an internal vertex
unless $G =  S_1 \cup S_2$.

We claim further (compare Lemma 4.17 of \cite{owe-pow} ) that each of the two
separation blocks $S_1(w)$ and $S_2(w)$ has at least two internal
vertices which are adjacent to the vertex $w$. For suppose to the
contrary that w has only one neighbour in $C_i(w)$. Then the proper
subgraph induced by the vertices of $S_i(w)\backslash w$ has freedom number
$f_i-1=0$ and $3$ vertices of attachment $\{w',x, y\}$. Therefore, it
cannot have an internal vertex so it is the $3$-cycle of edges
$(w'x), (w'y)$ and $(xy)$. This means the vertex $w'$ has degree $3$ and
is on a $3$-cycle containing $(xy)$ which was excluded at the start of
the proof.

Note that there may be several different vertices $w$ such that $(x,
y, w)$ separate $G$. Let $w$ and $w'$ be two vertices such that $(x, y, w)$
and $(x, y, w')$ both separate $G$ and suppose that $G$ contains the
edge $(ww')$. Then the vertex pair $\{x, y\}$ does not separate the
graph $G\backslash (ww')$ because this would imply that one component of the
separation set $(x,y,w)$ would have $w$ adjacent to only one internal
vertex $w'$.

Now consider a complete separation of $G/e$ at all separation pairs
as described in Lemma 5. By this lemma every separation block
is either a $3$-cycle or a $3$-connected graph. We will show that in
fact each of the separation blocks is a $3$-connected graph.

Since $G$ is $3$-connected every separation pair of $G/e$ consists of
the vertex $v$ which results from combining the vertices $x$ and $y$ of
$G$ and a vertex $w$ as described above. Thus every virtual edge is
incident to the vertex $v$. Since neither $(xw)$ nor $(yw)$ is an edge
of $G, (vw)$ is not an edge of $G/e$ and so every vertex in every
separation block which is a vertex of a separation pair of $G$
gets a virtual edge. Every vertex of $G/e$ has degree at least $3$ and
so any vertex in any separation block with degree less than
three is incident to a virtual edge. In particular there are no
$3$-cycles in which there are two original (non-virtual) edges at a
vertex. Since all virtual edges are adjacent to the vertex $ v$ the
only possible $3$-cycle consists of the vertex $v$, virtual edges $(vw)$
and $(vw')$ and an original edge $(ww')$. This means that the vertex $v$
separates $(G/e)\backslash (ww')$. This implies that $(x, y, w)$ and $(x,y,w')$
are distinct separation sets for $G$ and that the vertex pair $\{x,y\}$
separates $G\backslash (ww')$. This contradicts the remark above and shows
that all separation blocks of $G/e$ are $3$-connected graphs.
\hfill $\Box$

\medskip

These Lemmas now provide the basis for our proof of the main
theorem
\medskip

\noindent {\bf Theorem 1}: {\it Conjecture B implies Conjecture A.}
\medskip

\noindent {\it Proof.} Suppose to the contrary that Conjecture B is true and that
$G$ is a $3$-connected maximally independent graph which is not basic
and which is generically RS. Suppose further that $G$ is a vertex
minimal graph with this property. Then $|G| >6$ since the only
$3$-connected maximally independent graph which is not basic with
$|G| \le 6$ is the doublet and this is not RS.

Since $G$ is $3$-connected, every proper subgraph of $G$ has more
than two vertices of attachment in $G$. Let  $\{M\}$  be the set of
maximally independent proper subgraphs of $G$ and note that  $\{M\}$
is not empty since $G$ is not basic. If  $\{M\}$  contains a subgraph, $S$
which has an internal vertex then choose $R$ from  $\{M\}$
such that $R$ contains $S$ and $R$ is maximal. Otherwise simply
choose $R$ from  $\{M\}$  to
be maximal. In the former case, the internal vertex of $S$ is also
an internal vertex of $R$.

Let  $G\# R$ be the graph derived from $G$ as
described in Lemma 4. Then by Lemma 4  $G\# R$ is $3$-connected
and maximally independent and is generically RS. If $R$ has
any internal vertices, then $| G\# R| < |G|$ and this contradicts
the requirement that $G$ is minimal.

Otherwise $|G\# R| = |G| > 6$ and we may assume that $G$ and
thus  $G\# R$ has no maximally independent subgraphs which
contain internal vertices. By Lemma 6 either  $G\# R$ has a
contractible edge $f$ such that $( G\# R)/f$ is $3$-connected or
all of the separation blocks of $( G\# R)/e$ (including
all their virtual edges) are $3$-connected. In the first case
$( G\# R)/f$ is RS by Theorem 6.1 of  \cite{owe-pow}  which is quoted above
and $|( G\# R)/f | = |G|-1$. This contradicts the requirement
that $G$ is minimal. In the second case, we conclude from
Lemma 5 that at least one of the separation blocks $B_1$
of  $( G\# R)/e$ contains no redundant edges and so remains
$3$-connected when redundant edges are ignored. Then by Lemma
 5 $G$ is generically RS only if the $3$-connected graph $B_1$ is
 generically RS and since $|B_1| < |G|$ this contradicts the
 requirement that $G$ is minimal.
\hfill $\Box$

\section{K(3,3) is not Soluble by Radicals}

Let us label the vertices  of $K(3,3)$ as shown in figure 1 and
select  the coordinates of vertices $1$ and $2$ to be $(0,0)$ and $(0,1)$
respectively. Then the polynomials  represented by $K(3,3)$ are:
\[
x_3^2+y_3^2-d_1,~~~~ (x_4-1)^2+y_4^2-d_3,~~~~ x_5^2+y_5^2-d_2,
\]
\[(x_6-1)^2+y_6^2-d_4,~~~~
(x_3-x_4)^2+(y_3-y_4)^2-d_5,~~~~ (x_4-x_5)^2+(y_4-y_5)^2-d_6,
\]
and
\[
(x_4-x_6)^2+(y_4-y_6)^2-d_7,
~~~~
 (x_3-x_4)^2+(y_3-y_6)^2-d_8.
\]
For each choice of real algebraically independent  (squared)
distances  $d_1, \dots , d_8$, these equations determine a zero dimensional
complex affine variety.

In order to prove that the zeros of these equations are
non-radical for generic distances, we need to show that a
polynomial generator for the elimination ideal in one of the
variables is a polynomial whose Galois group is not a soluble
group. Although this can in principle be done by treating the
distances as generic parameters and performing the algebra in the
field $\mathbb{Q}(\{d\})$, in practice the Maple software package which we use
is unable to complete the calculation. We therefore choose
specialised rational values for the distances, perform the
calculations in $\mathbb{Q}$ and then use the specialisation Theorem 7.2 of
\cite{owe-pow}  to make a conclusion for generic distances. We select the
specialised values for the (squared) distances to be $d_1=d_2=d_3=d_4=1,
d_5=1/4,
d_6=4, d_7=9/16$ and $d_8= 9/4.$

The fifth equation and its three successors admit the squared form
\[
(d_5-(x_3-x_4)^2+y_3^2+y_4^2)^2-4y_3^2y_4^2=0,
\]
which in turn yields an equation in $x_3$ and $x_4$ alone on
substituting for $y_3^2$ and $y_4^2$ from the first two equations. In this
way  we obtain a system $\{g\} = \{g_1,g_2,g_3,g_4\}$ of four quartic
equations in $x_3,x_4,x_5,x_6$ and the squared distances. Computing the
resultants of pairs of these equations with respect to $x_4, x_5$ and
$x_6$ gives a polynomial for $x_3$ of degree $20$ for the specialised
distance values given. This polynomial is known to be a
specialisation of a polynomial with generic distances in the
elimination ideal for $x_3$ by Theorem 8.3 of \cite{owe-pow} . This polynomial
factors into three polynomials, namely (with $x = x_3)$, $(x-1)^6$ and
the following
polynomials
of degree 6 and degree 8;
\[
   87733791129600 x^6  - 280160493061120 x^5
      + 486601784497152 x^4 -
\]
\[
      581731370400244 x^3
         + 396516248769992 x^2 - 110509387701405 x + 1912924250825
\]
and \[
  19741148184576 x^8  - 103544588664832 x^7  + 29867097677824 x^6 +
\]
\[      440356364853504 x^5
       - 761674146310464 x^4
         + 517152016022904 x^3 -
\]
\[
215063281430796 x^2
         + 118596291789193 x - 45476733930709
\]
 The factor $(x-1)$  does not extend to give
a zero of $K(3,3)$ because this would require $x_3=1,y_3=0$ so that
points 2 and 3 are coincident which would require $d_3=d_5$. Maple
shows that the degree 6 and degree 8 polynomials have a
non-soluble Galois group. Then the specialisation Theorem 7.2 of
\cite{owe-pow}  allows us to conclude that the graph $K(3,3)$ is generically not
radically soluble.

\end{document}